\pageno=1
\magnification=1200
\vglue 1in
\hsize=11.9cm
\hoffset=0.9cm
\voffset=0.0 truecm  
\headline={\hss}
\footline={\tenrm\ifodd\pageno
\hss\folio\else
\folio\hss\fi}

\def \odd {{\rm odd}}
\def \even {{\rm even}}

\def\D{{\cal D}}
\def\L{{\cal L}}
\def\Nn{{\cal N}}

\def\O{{\cal O}}

\def\litem{\par\noindent\hangindent=\parindent}

\def\picture #1 by #2 (#3){
  \vbox to #2{
    \hrule width #1 height 0pt depth 0pt
    \vfill\special{picture #3}}} 
\def\scaledpicture #1by #2 (#3 scaled #4){{
  \dimen0=#1 \dimen1=#2
  \divide\dimen0 by 1000 \multiply\dimen0 by #4
  \divide\dimen1 by 1000 \multiply\dimen1 by #4
  \picture \dimen0 by \dimen1 (#3 scaled #4)}}

\def\O {{\cal O}}

\let\PP=\par
\let\mPP=\medbreak
\let\bPP=\bigbreak
\def\LP {\par\noindent}
\def\sLP {\smallskip\noindent}
\def\mLP {\medskip\noindent}
\def\bLP {\bigskip\noindent}
\def\br {\hfill\break}

\def\si {\sigma}

\def \odd {{\rm odd}}
\def \even {{\rm even}}

\def\part {\partial}

\def\Cov{{\rm Cov}}

\def\ltextindent#1{\hbox to \hangindent{#1\hss}\ignorespaces}

\lineskip=0pt

\def\footnote{  } {  }
\def\footnoterule{\kern-3pt
\hrule width\hsize\kern 2.6pt}
\font\eightrm=cmr8
\long\def\voetnoot#1#2{{\baselineskip=9pt
\setbox\strutbox=\hbox{\vrule height 7pt depth 2pt width 0pt}
\eightrm\footnote{#1}{#2}}}

\catcode`\@=12

\def\es{{1\over n-k+1}\sum\limits_{j=0}^{n-k}}
\def\est{{2\over n-k+1}\sum\limits_{j=0}^{n-k}}

\def\ets{{\lim\limits_{n\to\infty}}}
\def\etn{{E\bigl(c^TZ_n\bigr)}}
\def\Cov{{\rm Cov}}
\def\estoi{{2\over n-k+1}\sum\limits_{h=0}^{n-k}}

\centerline{ D\ I\ S\ C\ R\ E\ T\ E\ ~ S\ P\ A\ C\ I\ N\ G\ S}
\vskip1cm
{\obeylines \everypar {\hfil}\parindent=0pt 
Chris A.J.~Klaassen and J.~Theo Runnenburg
Korteweg - de Vries Institute for Mathematics
University of Amsterdam
The Netherlands}
\vskip2cm
\centerline{November 20, 2001}
\vskip2cm
\noindent
{\bf Summary}
\mPP
Consider a string of \ $n$ \ positions, i.e. a discrete string of
length \ $n$. Units of length \ $k$ \ are placed at random on this string in such
a way that they do not overlap, and as
often as possible, i.e. until all spacings between neighboring units have
length less than \ $k$. When centered and scaled by \ $n^{-1/2}$ \ the
resulting numbers of spacings of length \ $1, 2,\ldots ,k-1$ \ have
simultaneously a limiting normal distribution as \ $n\to\infty$.  This is
proved by the classical method of moments.
\vfill
{\parindent 0pt 
{{\sl AMS 1991 subject classifications.}} Primary 60F05 ; secondary 60D05\par}
{\parindent=4.5truecm
\item
\par}
\eject
\LP
{\bf 1. Introduction and main result}
\mPP
Imagine a row of \ $n$ \ hooks on which one hangs hats at random. Each hat
needs \ $k$ \ adjacent hooks.  Successive hats are put on the hooks by a random
device that  selects the next set of \ $k$ \ adjacent hooks to be covered by a
hat with equal probabilities from the remaining sets of \ $k$ \ adjacent
hooks.  This process continues until all sequences of adjacent free hooks have
size less than \ $k$.
\mPP
Let \ $X_{n,j}$ \ be the resulting number of sequences of adjacent free
hooks of exactly size \ $j$.  Necessarily \ $j$ \ ranges from \ $1$ \ to \ $k-1$. 
In the formulation of the problem given in the Summary, \ $X_{n,j}$ \ is the
number of spacings of length \ $j$. In this paper asymptotic normality of \
$X_n=\bigl(X_{n,1},\ldots ,X_{n,k-1}\bigr)^T$ \ is proved.
\mLP
{\bf Theorem 1.1} 
{\sl Let \ $X_n$ \ be the vector of numbers of spacings
of length \ $j$, \hfil\break  $j=1,\ldots, k-1$. Then
$$Z_n=n^{-1/2}\bigl(X_n-EX_n\bigr){\buildrel\D\over\longrightarrow\Nn}(0,\Sigma)
\qquad  ,{\rm as} \enskip n\to\infty,\leqno (1.1)$$ 
and
$$n^{-1/2}\bigl(EX_n-n\theta\bigr) \longrightarrow 0
\qquad  ,{\rm as} \enskip n\to\infty,\leqno (1.2)$$ 
hold. The components \ $\theta_j, ~ j=1,\ldots, k-1,$ of $\theta$ are described
in {\rm Lemma 3.1} and \ $\Nn(0,\Sigma)$ \ 
denotes the centered normal distribution with covariance
matrix \ $\Sigma =\bigl(\sigma_{ij}\bigr)^{k-1}_{i,j=1}$ \ and with \ $\sigma_{ij}$ \
as defined in {\rm (3.13)} through {\rm (3.15)}.} 
\mPP
 An outline of our proof of this limit theorem is given in Section 2. The
limit behavior of expectation vector and covariance matrix of \ $X_n$ \ is studied
in Section 3.
The case \ $k=2$ \ has been treated along different lines by Runnenburg (1982),
who also gives references to literature about this problem. As Runnenburg (1982)
mentions, the random device we study here is
relevant to problems in chemistry. Mackenzie (1962)
describes its application to complete adsorption of linear 
molecules into parallel troughs (valleys) on suitable crystal surfaces.
Here each trough is represented by our row of hooks and
our hats symbolize the linear molecules. Mackenzie (1962)
has derived the asymptotic behavior of the mean $EX_n$ by a
similar technique as in our Lemma 3.1. The strength of our
approach is in Lemma 2.1, which enables us to derive the 
asymptotic normality from the convergence of the first two moments.
Our results are compared to these results from literature in Section 4.
\mPP
 A continuous version of this problem has been suggested by R\'{e}nyi (1958). 
It is called the parking problem and it has been studied by
Dvoretzky and Robbins (1964).
Coffman, Flatto, and Jelenkovi\'c (2000) study the number
of vacant intervals in this parking problem when the cars
arrive according to a Poisson process in time and place. 
They consider the asymptotic behavior of mean and variance 
of this number of intervals at fixed time as the available 
space (the length of the row of hooks in the discretized version)
grows beyond bounds.
\eject

\bLP
{\bf 2. Heuristics and a key lemma}
\mPP
Clearly, when placing the first hat there are \ $n-k+1$ \ possibilities,
which all have the same probability \ $1/(n-k+1)$. If the first hat occupies
positions \ $j+1$ \ through \ $j+k$ \ then two rows of hooks are left, one of
size \ $j$ \ and the other one of size \ $n-k-j$. 
Let \ $J$ \ be a random
variable uniformly distributed on \ $\{0,1,\ldots ,n-k\}$ \ and let the \
$(k-1)$-vector \ $\tilde X_i$ \ have the same distribution as \ $X_i$. 
We assume
that \ $J,~ X_0,X_1,\ldots,X_n,  {\tilde X}_0,{\tilde X}_1,\ldots,{\tilde
X}_n$ \ are independent. 
The above reflection shows (cf. Lemma 2.4 of Runnenburg
(1982)) 
$$\L(X_n)=\L\bigl(X_J+{\tilde X}_{n-k-J}\bigr),\qquad n=k, \ k+1,\ldots .\leqno
(2.1)$$
\mPP
In terms of characteristic functions this can be formulated simpler:
$$\eqalign{Ee^{is^TX_n}=&{1\over
n-k+1}\sum\limits_{j=0}^{n-k}\biggl(Ee^{is^TX_j}\biggr)
\biggl(Ee^{is^TX_{n-k-j}}\biggr),\cr   &s\in {\bf R}^{k-1}, \quad n=k, k+1,\dots\,
.\cr}\leqno (2.2)$$
With \ $s=n^{-1/2}t$ \ and \ $Z_n$ \ as in the theorem this becomes
$$\eqalign{ Ee^{it^TZ_n}&={1\over
n-k+1}\sum_{j=0}^{n-k}\bigl(Ee^{it^T(j/n)^{1/2}Z_j}\bigr)\cr
&\bigl(Ee^{it^T((n-k-j)/n)^{1/2}Z_{n-k-j}}\bigr)
e^{in^{-1/2}t^T{(EX_j+EX_{n-k-j}-EX_n)}}.\cr}\leqno(2.3)$$ 
\noindent Here \ $Z_0$ \ is degenerate at 0. 
\mPP
If we could prove that for some finite constant \ $B$ 
$$\|EX_j+EX_{n-k-j}-EX_n\|\le B, \quad j = 0,\ldots ,n-k,\quad  n=k,k+1,\ldots
,\leqno (2.4)$$
holds with \ $\|\cdot\|$ \ the Euclidean norm and that \ $\{Z_n\}$ \ has a limit
distribution with characteristic function \ $\psi$ \ then we could deduce from (2.3)
$$\psi(t)=\int\limits_0^1\psi({\sqrt u}\,t)\psi({\sqrt {1-u}}\,t)du, \quad
t\in{\bf R}^{k-1}.\leqno (2.5)$$
\mPP
 By Lemma A.1 from the Appendix this integral equation implies that $\{Z_n\}$
is asymptotically normal, but the asymptotic covariance structure
is not determined in this way, not even if we would know the limit behavior 
of the covariance matrix of $Z_n$. Therefore,
we will not try to elaborate the details of this approach, but we
will use the classical method of moments instead,
which does yield the asymptotic covariance matrix.
\mPP
 This method of moments is based on the theorem of Fr\'echet
and Shohat (1931) about the celebrated moment convergence problem,
which states that if all moments of $\{Z_n\}$ converge to the
corresponding moments of $Z$ and if the moments of $Z$ determine
its distribution uniquely, then $\{Z_n\}$ converges to $Z$ in
distribution; see 11.4.C, p.187, of Lo\`eve (1977) or
Tak\'acs (1991).
Our key-stone for the method of moments is the fact that 
the normal distribution is indeed determined uniquely by
its moments. By considering $\{c^T Z_n\}$ and $c^T Z$ 
for all $c$ we see that the multidimensional result is obtained
from the one-dimensional result.
\mPP
In Section 3 we shall show the validity of (2.4).
Let ~${\tilde Z}_n \br =n^{-1/2}({\tilde X}_n-E{\tilde X}_n), n=1,2,\dots, {\tilde
Z}_0=0$~ a.s.  Note that it follows from (2.3) and (2.4) that for all \ $c\in {\bf R}^{k-1}$ \ and
all \ $m\in{\bf N}$, \ once we have \ $\sup_j \ E|c^TZ_j|^{m-1}<\infty$, \ 
then we obtain 
$$\eqalign{E\bigl(&c^TZ_n\bigr)^m\cr 
&=\es E\biggl({\sqrt{j\over
n}}\enskip c^T Z_j+{\sqrt{n-k-j\over n}} \ c^T{\tilde Z}_{n-k-j}\cr
&\qquad\qquad\qquad\qquad+{1\over{\sqrt n}} \ 
c^T\bigl(EX_j+EX_{n-k-j}-EX_n\bigr)\biggr)^m\cr  
&=\es E\biggl\{\sum\limits_{i=0}^m{m\atopwithdelims () i}\biggl({\sqrt{j\over
n}} \ c^TZ_j\biggr)^i\cr
&\qquad\qquad\biggl({\sqrt {n-k-j\over n}} \
c^T{\tilde Z}_{n-k-j}\biggr)^{m-i}\biggr\}+\O\biggl(n^{-1/2}\biggr)\cr
&=\sum\limits_{i=1}^{m-1}{m\atopwithdelims () i}\es \biggl({j\over
n}\biggr)^{i/2}\biggl({n-k-j\over n}\biggr)^{(m-i)/2}
E\biggl(c^TZ_j\biggr)^i\cr
&\qquad\cdot E\biggl(c^TZ_{n-k-j}\biggr)^{m-i}+\est\biggl({j\over n}\biggr)^{m/2}
E\biggl(c^TZ_j\biggr)^m\cr 
&\qquad +\O\biggl(n^{-1/2}\biggr), \ {\rm as} \
n\to\infty,\cr}\leqno (2.6)$$
and similarly
$$\eqalign{&E|c^TZ_n|^m\cr
&\quad\le \es E\biggl({\sqrt{j\over
n}}\big|c^T Z_j\big|+{\sqrt{n-k-j\over n}}|c^T{\tilde Z}_{n-k-j}|\cr
&\qquad\qquad\qquad\qquad\qquad+{1\over{\sqrt n}} \ 
\big|c^T\bigl(EX_j+EX_{n-k-j}-EX_n\bigr)\big|\biggr)^m\cr  
&=\sum\limits_{i=1}^{m-1}{m\atopwithdelims () i}\es \bigg({j\over n}\bigg)^{i/2}
\bigg({n-k-j\over n}\bigg)^{(m-i)/2} E|c^TZ_j|^i\cr
&\quad\cdot E\big|c^TZ_{n-k-j}\big|^{m-i}+\est\biggl({j\over n}\biggr)^{m/2}
E\big|c^TZ_j\big|^m\cr 
&\quad +\O\big(n^{-1/2}\big). \cr}\leqno (2.7)$$
These relations between the (absolute) moments of \ $c^TZ_n$ \ enable us to derive the
asymptotic behavior of all moments of order \ $m>2$ \ from that of the
moments of orders \ $1$ \ and \ $2$, \ by the following key result.
\mLP
{\bf Lemma 2.1} {\sl Let \ $\{a_n\}^\infty_{n=0}$ \ be a sequence of reals. If for some \
$k\in{\bf N}, \ \alpha\in{\bf R},$ and \ $\beta>1$
$$a_n=\alpha +{2\over n-k+1}\sum\limits_{j=0}^{n-k}\biggl({j\over
n}\biggr)^\beta a_j+{\scriptscriptstyle{\cal O}}(1), \ {\rm as} \ n\to\infty, 
\leqno (2.8)$$
holds, then \ $\{a_n\}$ \ converges and
$$\lim_{n\to\infty} a_n = \alpha\, {\beta+1\over\beta-1}\, \leqno (2.9)$$
holds. Moreover, if the elements ~$a_n$~ of a possibly other sequence
\ $\{a_n\}^\infty_{n=0}$ \ are nonnegative and bounded from above by the right-hand side of {\rm (2.8)} for some \
$k\in{\bf N}, \ \alpha\in{\bf R},$ and \ $\beta>1$, then ~$\sup_na_n$~ is finite.}    
\sLP
{\bf Proof.} If \ ${\overline {\lim}}_{n\to\infty} |a_n|=\infty$, \ then there
exists a sequence \ $\{n_i\}, \ n_i\to\infty$ \ as \ $i\to\infty$, \ such that
$$\eqalign{&|a_{n_i}|=\sup\{|a_h|:h\le n_i\},\cr
&|a_{n_i}|\to\infty.\cr}$$
\mPP
By (2.8) this yields
$$\eqalign{1=\biggl|{a_{n_i}\over
a_{n_i}}\biggr|&\le\biggl|{2\over n_i-k+1}\sum\limits_{j=0}^{n_i-k}
\biggl({j\over n_i}\biggr)^\beta{a_j\over
|a_{n_i}|}\biggr|+{\scriptscriptstyle{\cal O}}(1)\cr
&\le {2\over n_i-k+1}\sum\limits_{j=0}^{n_i-k}\biggl({j\over n_i}\biggr)^
\beta+{\scriptscriptstyle{\cal O}}(1)=2\int\limits_0^1 x^\beta dx+{\scriptscriptstyle{\cal O}}(1)\cr
&={2\over\beta +1}+{\scriptscriptstyle{\cal O}}(1),\cr}$$
which is a contradiction in view of \ $\beta>1$. 
Consequently, we have
$$ {\mathop{\overline{\lim }}\limits_{n\to \infty }}|a_n|<\infty, \ {\rm and} \
\sup\limits_n|a_n|=a<\infty.\leqno (2.10)$$
Incidentally, this proves the second statement of the Lemma.
Let \ $m_n$ \ be a sequence of integers tending slowly to infinity; more
precisely, \ $m_n\to\infty$ \ and \ $m_n/n\to 0$ \ as \ $n\to\infty$. 
By (2.8) and (2.10) we obtain
$$\eqalign{{\mathop{\overline{\lim }}\limits_{n\to \infty }} a_n \le 
&\alpha + {\mathop{\overline{\lim }}\limits_{n\to \infty }}\biggl\{ 
{2\over{n-k+1}}\sum\limits_{j=m_n}^{n-k}\biggl({j\over
n}\biggr)^\beta\sup_{h\ge m_n} a_h+2a{{m_n}\over{n-k+1}}\biggr\}\cr  
=&\alpha +{\mathop{\overline{\lim}}\limits_{n\to \infty
}}\biggl\{2\int\limits_0^1 x^\beta dx \sup_{h\ge m_n} a_h\biggr\}\cr  
=&\alpha+{2\over\beta +1}{\mathop{\overline{\lim }}\limits_{n\to \infty }} 
a_n.\cr}\leqno(2.11)$$ 
 Consequently, we have
$${\mathop{\overline{\lim }}\limits_{n\to \infty }}\, a_n\le\alpha\,{\beta
+1\over\beta -1}.\leqno (2.12)$$
Similarly, we obtain
$$\eqalign{{\mathop{\underline{\lim }}\limits_{n\to \infty }}\, a_n&\ge\alpha +
{\mathop{\underline{\lim }}\limits_{n\to \infty }}\biggl\{{2\over\beta
+1}\inf_{h\ge m_n} a_h-2a\,{m_n\over n-k+1}\biggr\}\cr
&=\alpha +{2\over\beta +1} {\mathop{\underline{\lim }}\limits_{n\to \infty
}}\, a_n,\cr}$$
and hence
$${\mathop{\underline{\lim }}\limits_{n\to \infty }}\, a_n\ge\alpha\, {\beta
+1\over\beta -1}.\leqno(2.13)$$
\hfill QED
\bPP 
 From (2.6) and Lemma 2.1 with \ $\beta=m/2>1$ \ it follows
that all moments converge once the first two (absolute) moments do and more
explicitly that 
$$\eqalign{\lim\limits_{n\to\infty} E\biggl(c^TZ_n\biggr)^m
=&{m+2\over m-2}\sum_{i=1}^{m-1}{m\choose i}\int\limits_0^1x^{i/2}(1-x)^{(m-i)/2}dx\cr 
&\cdot\lim\limits_{n\to\infty}E\biggl(c^TZ_n\biggr)^i
\lim\limits_{n\to\infty}E\biggl(c^TZ_n\biggr)^{m-i}\cr}\leqno(2.14)$$
holds. Note that existence (and boundedness in $n$) of the odd moments 
is verified here by application of the second part of Lemma 2.1 to the corresponding
absolute moments satisfying (2.7). Of course, the first moment vanishes. Assume
$$\ets\,E(c^TZ_n)^2=\sigma^2<\infty\leqno(2.15)$$
and note that this implies the existence and boundedness in $n$ of the first
absolute moments. Denote \ $\ets\,\etn^m$ \ by \ $\mu_m$ \ for \ $m>2$. 
 From (2.14) and by induction on \ $m$ \ we arrive at
$$\eqalign{\mu_m =\cases{0  &\quad $\odd$\cr
\cr
\qquad\qquad\qquad $if \ $m$ \ is$  \cr
\cr
2^{-m/2}m!\sigma^m/(m/2)!   &\quad $\even$,\cr}}\leqno(2.16)$$
which are the moments of a normal distribution with mean 0 and variance $\sigma^2$.
This yields asymptotic normality of \ $c^TZ_n$ \ by the method of moments as
discussed in Section 1, \ since the normal distribution is
determined by its moments; see e.g. Feller (1971), Examples VIII.6(b) and
VIII.1(e), pages 269 and 251, 252.
Consequently, the proof of (1.1) is complete once we have proved (2.4)
and (2.15) with \ $\sigma^2=c^T\Sigma c$. This will be done in the next section.

\bLP
{\bf 3. Asymptotics for the first two moments}
\mPP 
Asymptotically \ $EX_n$ \ behaves as follows.
\mLP
{\bf Lemma 3.1} {\sl With the notation  
$$e_k(y)=\exp\bigl\{2[y+y^2/2+\ldots +y^{k-1}/(k-1)]\bigr\}$$ 
and
$$\theta_j=2\bigl(e_k(1)\bigr)^{-1}\int\limits_0^1 (1-y)y^je_k(y)dy, \quad
j=1,\dots,k-1,$$  
there exists for every \ $R>1$ \ a constant \ $B_R$ \ with}
$$|EX_{n,j}-(n+k)\theta_j|\le B_R\,R^{-n},\quad j=1,\dots,k-1,\quad  n=1,2,\ldots
\, .\leqno(3.1)$$ 
\PP 
Note that (3.1) implies (2.4) and (1.2), and enables us to conclude from (1.1) 
that also
$$n^{-1/2}\biggl(X_n-n\theta\biggr){\buildrel\D\over\longrightarrow\Nn}
(0,\Sigma)\leqno(3.2)$$ 
holds. Therefore, the proof of Theorem 1.1 is complete once (1.1) has been shown.
\mLP
{\bf Proof of Lemma 3.1}
Fix \ $c\in{\bf R}^{k-1}$ \ and let \ $\gamma_n=Ec^TX_n$.  By
(2.1) or (2.2)
$$\eqalign{\gamma_n&=\es(\gamma_j+\gamma_{n-k-j})\cr
&=\est\gamma_j, \quad n=k, k+1,\ldots ,\cr}\leqno(3.3)$$
and hence
$$(n-k+1)\gamma_n=(n-k)\gamma_{n-1}+2\gamma_{n-k}, \quad n=k+1, k+2,\dots
,\leqno(3.4)$$
hold. 
Using generating functions similar to Overdijk (1981)
(see also Mackenzie (1962)) we define 
$$G(z)=\sum_{n=k+1}^\infty \ \gamma_nz^{n-k+1} \leqno(3.5)$$ 
for \ $|z|<1$ \ (note $\gamma_n= {\cal O}(n)$). 
Multiplying (3.4) by \ $z^{n-k}$ \ and summing over
\ $n=k+1, k+2,\ldots$ \ we obtain
$$(1-z)G'(z)=2z^{k-1}G(z)+2[\gamma_1z+\ldots+\gamma_kz^k]+\gamma_kz.$$
In fact, \ $\gamma_k=0$ \ since \ $X_k=0$ \ a.s., and \ $\gamma_i=c_i$ \ for \
$i=1,\ldots ,k-1$. Solving this differential equation under the side condition \
$G(0)=0$ \ we obtain
$$\eqalign{G(z)&=(1-z)^{-2}\psi(z), \cr
\psi(z)&=2\bigl(e_k(z))^{-1}\int\limits_0^z \bigl[c_1y+\ldots
+c_{k-1}y^{k-1}\bigr](1-y)e_k(y)dy.\cr}\leqno(3.6)$$
Note that \ $\psi(z),~z\in{\bf C}$, \ is an entire function and
write \ $\psi(z)=\Sigma_{i=2}^\infty a_iz^i$. 
Then (3.6) yields
$$G(z)=\sum\limits_{n=k+1}^\infty\,\sum\limits_{i=2}^{n-k+1}(n-k+2-i)
a_iz^{n-k+1}\leqno(3.7)$$
and hence by comparison to (3.5), for ~$n=k+1, k+2,\dots\,$,
$$\eqalign{\gamma_n&=\sum\limits_{i=2}^{n-k+1}(n-k+2-i)a_i\cr
&=\sum\limits_{i=2}^\infty(n-k+2-i)a_i-\sum\limits_{i=n-k+3}^\infty(n-k+2-i)a_i\cr
&=(n-k+2)\psi(1)-\psi'(1)+R_n\cr}\leqno(3.8)$$
with
$$\bigl|R_n\bigr|=\big|\sum\limits_{i=n-k+3}^
\infty\bigl(i-(n-k+2))a_i\big|\le \sum\limits_{i=n-k+3}^\infty
i\big|a_i\big|.\leqno(3.9)$$ 
Since \ $\psi(z),~z\in{\bf C}$, \ is an entire function, 
\ $z\psi'(z)=\sum_{i=2}^\infty ia_iz^i$ \ is entire as well and by Cauchy's
root test for the radius of convergence of analytic functions we obtain
$${\mathop{\overline{\lim }}\limits_{i\to \infty }}|ia_i|^{1/i}=0.\leqno(3.10)$$
Consequently, for arbitrary \ $R>1$ \ and \ $n$ \ sufficiently large we have
$$\sum_{i=n}^\infty\, i|a_i|\le\sum_{i=n}^\infty\,R^{-i}={R^{-n}\over
1-R^{-1}}.\leqno(3.11)$$
Finally, note
$$(n-k+2)\psi(1)-\psi'(1)=(n+k)\psi(1)=(n+k)c^T\theta,\leqno(3.12)$$
which by (3.8),(3.9), and (3.11) implies the lemma. 
\hfill QED
\bPP 
For our study of the (asymptotic) covariance structure of \ $X_n$ \ we need
some extra notation. Let
$$G_i(z)=2(1-z)^{-2}\bigl(e_k(z)\bigr)^{-1}
\int\limits_0^z y^i(1-y)e_k(y)dy\leqno(3.13)$$
be the generating function as in (3.5) through (3.6) of the means of
$X_{n,i} \,,i=1,\ldots,k-1, n=k+1,k+2,\ldots,$ and define
$$\eqalign{H_{ij}(z)= & (1-z)z^i{\bf 1}_{[i=j]} +(1-z)^2(z^i+z^{k-1}G_i(z))
(z^j+z^{k-1}G_j(z))\cr
&-\theta_i\theta_j (1-z)^{-2}\biggl(3+(4k-5)(1-z)+2(k-1)^2(1-z)^2 \cr
& -2k^2(1-z)^4 -\bigl(2+(4k-3)(1-z)+(2k-1)^2(1-z)^2 \cr
& -4k^2(1-z)^3\bigr)z^k\biggr), \qquad i,j=1,\ldots ,k-1.} \leqno(3.14)$$
We will show that
$$\sigma_{ij}=2\bigl(e_k(1)\bigr)^{-1}\int\limits_0^1 H_{ij}(y)e_k(y)dy
\leqno(3.15)$$
is well defined and even that
\mLP
{\bf Lemma 3.2} {\sl For \ $i,j=1,\ldots,k-1$, \ and for every \ $R>1$,
$$\Cov(X_{n,i},X_{n,j})=(n+k)\sigma_{ij}+{\cal O}(R^{-n}), \
{\rm as} \ n\to\infty, \leqno(3.16)$$  
holds and consequently \ $\Sigma$ \ of Theorem 1.1 satisfies}  
$$\Sigma=\bigl(\sigma_{ij}\bigr)_{i,j=1}^{k-1}.\leqno(3.17)$$ 
\mLP
{\bf Proof.} Fix \ $i,j\in\{1,\ldots ,k-1\}$ \ and let \
$c_n=EX_{n,i}X_{n,j}-(n+k)^2\theta_i\theta_j$\ for \ $n=0,1,\dots$. 
By (2.1) or (2.2) we obtain
$$c_n=\estoi
\{c_h+(h+k)^2\theta_i\theta_j+EX_{h,i}EX_{n-k-h,j}\}
-(n+k)^2\theta_i\theta_j\leqno(3.18)$$
and hence
$$(n-k+1)c_n=(n-k)c_{n-1}+2c_{n-k}+2\delta_{n,i,j}, \qquad n=k+1,\, k+2,\ldots
,\leqno(3.19)$$ with 
$$\eqalign{\delta_{n,i,j}= & \sum_{h=0}^{n-k-1}\{\epsilon_{h,i}\theta_j
+\theta_i\epsilon_{h,j}+\epsilon_{h,i}(\epsilon_{n-k-h,j}-\epsilon_{n-k-h-1,j})\} \cr
& +k\theta_i\epsilon_{n-k,j},}$$
$$\epsilon_{h,i}=EX_{h,i}-(h+k)\theta_i, \qquad i=1,\ldots,k-1,h=0,1,\ldots.$$
Tedious computation shows
$$\eqalign{\sum_{n=k+1}^{\infty}\delta_{n,i,j}z^{n-k} = 
& \sum_{h=0}^{\infty}\sum_{n=k+h+1}^{\infty} \{
\epsilon_{h,i}z^h\theta_jz^{n-k-h} +\theta_i\epsilon_{h,j}z^hz^{n-k-h} \cr
& + \epsilon_{h,i}z^h\epsilon_{n-k-h,j}z^{n-k-h} 
- \epsilon_{h,i}z^h\epsilon_{n-k-h-1,j}z^{n-k-h-1}z \} \cr
& +k\theta_i\sum_{h=1}^{\infty}\epsilon_{h,j}z^h \cr
= &(1-z)(z^i+z^{k-1}G_i(z))(z^j+z^{k-1}G_j(z)) \cr
& +\theta_i\theta_j\{k^2-(1-z)^{-3}(1+(k-1)(1-z))^2\}.}$$
By an argument similar to that to derive (3.6) through (3.8) from (3.4) we obtain from (3.19)
$$G_{ij}(z)=\sum\limits_{n=k+1}^\infty c_nz^{n-k+1}=(1-z)^{-2}\psi_{ij}(z)
\leqno(3.20)$$
with
$$\psi_{ij}(z)=2\bigl(e_k(z)\bigr)^{-1}\int\limits_0^z H_{ij}(y)e_k(y)dy=
\sum\limits_{h=3}^\infty b_hz^h, \leqno(3.21)$$
and consequently, as in (3.6) through (3.9) and (3.12), 
$$\eqalign{c_n&=(n-k+2)\psi_{ij}(1)-\psi'_{ij}(1)+R_n\cr
&=(n+k)\si_{ij}+R_n\cr}\leqno (3.22)$$
with \ $\psi_{ij}(1)=\si_{ij}$ \ given in (3.15) and with
$$|R_n|\le \sum\limits_{h=n-k+3}^\infty h|b_h|.\leqno (3.23)$$
To complete the argument along the lines of the end of the proof of Lemma 3.1, i.e.
along the lines of (3.10) we have to show that $\psi_{ij}(z), ~z\in{\bf C}$,
is entire. Indeed, it may be verified that, as $z \to 1$, the function 
$G_i(z), ~z\in{\bf C}$, behaves like 
$\theta_i(1-z)^{-2}(1+2(k-1)(1-z)) + {\cal O}(1)$
and hence $H_{ij}(z)$ like ${\cal O}(1)$. Consequently, 
$H_{ij}(z), ~z\in {\bf C}$, is entire, since $(1-z)^2G_i(z)$ is.
Together with Lemma 3.1 this implies that (3.22) yields (3.16).
\hfill QED
\mPP 
Note that by the argument of Section 2, (3.16) implies that all moments of
\ $c^TZ_n$ \ with \ $Z_n$ \ as in (1.1) converge to those of a normal
distribution. Consequently, Theorem 1.1 holds with \ $\Sigma$ \ as in
(3.17).

\bLP
{\bf 4. Comparison to literature}
\mPP
For $k=2$ and $i=j=1$, Lemma 3.1 holds with
$$\theta_1=2e^{-2}\int_0^1(1-y)ye^{2y}dy=e^{-2}.\leqno (4.1)$$
Furthermore, lengthy computations show
$$G_1(z)=(1-z)^{-2}e^{2z}-1,\leqno (4.2)$$
$$\eqalign{H_{11}(z)= & z(1-z)+(1-2(1-z)^{-1}+(1-z)^{-2})e^{-4z} \cr
&-e^{-4}\bigl((1-z)^{-2}+2(1-z)^{-1}+1+29(1-z) \cr
&-49(1-z)^2+16(1-z)^3\bigr),}\leqno(4.3)$$
and finally
$$\sigma_{11}=2e^{-2}\int_0^1 H_{11}(y)e^{2y}dy=4e^{-4}.\leqno(4.4)$$
These results agree with the expressions for the asymptotic mean and
variance obtained in (3.20) and (3.27) respectively, of Runnenburg (1982).
 Mackenzie (1962) has studied asymptotic mean and variance of the 
total vacant length 
$$V_n=\sum_{i=1}^{k-1}jX_{n,j}.\leqno (4.5)$$
Our Theorem 1.1 yields asymptotic normality of the standardized $V_n$.
Moreover, Lemma 3.1 and
$$\eqalign{e_k(1)\sum_{j=1}^{k-1}j\theta_j 
& = 2\int\limits_0^1 (1-y)(y+2y^2+\dots +(k-1)y^{k-1})e_k(y)dy \cr
& = 2\int\limits_0^1 \bigl((ky-(k-1)y^2)(1+y+\dots
+y^{k-2})-(k-1)y\bigr)e_k(y)dy \cr
& = [(ky-(k-1)y^2)e_k(y)]_0^1 -k\int_0^1e_k(y)dy \cr
& = e_k(1)- k\int\limits_0^1 e_k(y)dy } \leqno(4.6)$$
yield the existence for every $R>1$ of a constant $B_R$  with
$$|EV_n-(n+k)\bigl(1- k(e_k(1))^{-1}\int\limits_0^1 e_k(y)dy|\le B_R\,R^{-n},
n=1,2,\ldots\, .\leqno(4.7)$$ 
This result is in line with (A28) and (A25) of Mackenzie (1962).
An asymptotic expression with an exponentially small error as above can
be obtained also for the variance of $V_n$ via Lemma 3.2.
\bLP
{\bf Appendix}
\mPP
{\bf Lemma A.1}
{\sl Any random variable on ${\bf R}^d$ with characteristic function $\psi$
satisfying integral equation {\rm (2.5)}, i.e.
$$\psi(t)=\int\limits_0^1\psi({\sqrt u}\,t)\psi({\sqrt {1-u}}\,t)du, \quad
t\in{\bf R}^{k-1},$$
has a normal distribution.}
\mLP
{\bf Proof.}
Fix $a \in {\bf R}^d$ and define
$$ \chi_a(s) = \int_0^{\infty} e^{-sw} \psi(\sqrt w \,a) dw\,, s>0.\leqno (A.1)$$
By the boundedness of $\psi$ this Laplace transform is well defined and
by dominated convergence we obtain
$$ \chi^{'}_a(s) ={d \over ds} \chi_a(s)
   = -\int_0^{\infty} e^{-sw} w \psi(\sqrt w \,a) dw\,.\leqno (A.2)$$
Consequently, (2.5) yields
$$ \eqalign {-\chi^{'}_a(s) = &\int_0^{\infty} e^{-sw} w\int_0^1 \psi(\sqrt{uw}\,a)
   \psi(\sqrt{(1-u)w}\,a)du\,dw\cr
  = &\int_0^{\infty}\int_0^w e^{-sw} \psi(\sqrt v\,a)
   \psi(\sqrt{w-v}\,a)dv\,dw\cr
  = &\int_0^{\infty} \int_v^{\infty} e^{-s(w-v)} \psi(\sqrt{w-v}\,a)\,dw
    e^{-sv} \psi(\sqrt v \,a)\,dv\cr
  = &\chi_a^2(s)\cr}\leqno (A.3)$$
and hence there exists $b_a \in {\bf C}$ with
$$ \chi_a(s)={1 \over b_a + s}\,, s>0.$$
By Laplace inversion (cf. e.g. VII.6.6, p.233, of Feller (1971))
this implies
$$  \psi(\sqrt w\,a) = e^{-b_aw}\, ,w>0.$$
Consequently, a random $d$-vector $Z$ with characteristic
function $\psi$ satisfies
$$ Ee^{isa^TZ} = \psi(sa) = e^{-b_as^2}\, ,s>0,$$
and hence $a^TZ$ is normal for any $a \in {\bf R}^d$ and $Z$
is multivariate normal indeed. \hfil \hfill QED
\mPP

{\bf Acknowledgements}
We would like to thank Jaap Korevaar and Jan Wiegerinck
for suggesting the use of Laplace transforms in analyzing 
equation (2.5) and Jordan Stoyanov for the reference to
Tak\'acs (1991).
\mPP

{\bf  References}
\mLP
{\parindent=2cm
\litem{Coffman, E.G. Jr., L. Flatto, and P. Jelenkovi\'c} (2000).
Interval packing: the vacant interval distribution.
{\sl Ann. Appl. Probab.} {\bf 10}, 240--257. \par}
\sLP
{\parindent=2cm
\litem{Dvoretzky, A. and H. Robbins} (1964).
On the `parking' problem.
{\sl MTA Mat. Kut. Int. K\"{z}l.} {\bf 9}, 209--225. \par}
\sLP
{\parindent=2cm
\litem{Feller,~W.}(1971). 
{\sl An Introduction to Probability Theory and Its
Applications}, Vol.II, 2nd ed., Wiley, New York. \par}
\sLP
{\parindent=2cm
\litem{Fr\'echet, M. and J. Shohat} (1931).
A proof of the generalized second-limit theorem in the
theory of probability.
{\sl Trans. Amer. Math. Soc.} {\bf 33}, 533--543. \par}
\sLP
{\parindent=2cm
\litem{Lo\`eve, M.} (1977).
{\sl Probability Theory}, Vol.I, 4th ed., Springer, New York. \par}
\sLP
{\parindent=2cm
\litem{Overdijk, D.} (1981). 
Solution of Problem 102(101) of Runnenburg and Steutel (in
Dutch).  {\sl Statist.  Neerlandica} {\bf 35}, 175--178. \par}
\sLP
{\parindent=2cm
\litem{R\'{e}nyi, A.} (1958). 
On a one-dimensional random space-filling problem.
{\sl MTA Mat. Kut. Int. K\"{z}l.} {\bf 3}, 109--127. \par}
\sLP
{\parindent=2cm
\litem{Runnenburg,}\enskip J.Th. (1982). 
Asymptotic normality in vacancies on a line.
{\sl Statist. Neerlandica} {\bf 36}, 135--148. \par}
\sLP
{\parindent=2cm
\litem{Tak\'acs, L.} (1991).
A moment convergence theorem.
{\sl Amer. Math. Monthly} {\bf 98}, 742--746. \par}

\end